\documentclass[12pt]{amsart}
\usepackage{amscd}

\newtheorem{pr}{Proposition}
\newtheorem{lm}{Lemma}
\newtheorem{tm}{Theorem}
\newtheorem{cor}{Corollary}

\newcommand{\com}{{\mathbb C}}

\newcommand{\rarr}{\rightarrow}
\newcommand{\Q}{{\mathbb{Q}}}
\newcommand{\bbS}{{\mathbb{S}}}
\newcommand{\CC}{{\mathcal C}}
\newcommand{\scup}{\cup}
\newcommand{\eqq}{\stackrel{\sim}{=}}
\newcommand{\bpf}{\noindent {\em Proof.} }
\newcommand{\epf}{\qed \vspace{+10pt}}
\newcommand{\M}{{\overline M}}
\newcommand{\Spec}{{\rm Spec}}

\begin{document}
\title{Constructions of nontautological classes on
moduli spaces of curves}
\author{T. Graber and R. Pandharipande}
\date{10 March 2001}
\maketitle

\pagestyle{plain}
\baselineskip=16pt 

\setcounter{section}{-1}
\section{Introduction}
The {\em tautological rings} $R^*(\overline{M}_{g,n})$ are natural
subrings of the Chow rings of the Deligne-Mumford
moduli spaces of pointed curves:
\begin{equation}
\label{tautsys}
R^*(\overline{M}_{g,n}) \subset A^*(\overline{M}_{g,n})
\end{equation}
(the Chow rings are taken with $\mathbb{Q}$-coefficients).
The system of tautological subrings (\ref{tautsys}) is defined
to be the set of smallest $\Q$-subalgebras  satisfying the
following three properties [FP]:
\begin{enumerate}
\item[(i)] $R^*(\overline{M}_{g,n})$ contains
the cotangent line classes 
$$\psi_1, \ldots, \psi_n \in A^1(\overline{M}_{g,n}).$$
\item[(ii)] The system is closed under push-forward via
all maps forgetting markings:
$$\pi_*: R^*(\overline{M}_{g,n}) \rarr R^*(\overline{M}_{g,n-1}).$$
\item[(iii)] The system is closed under push-forward via
all gluing maps:
$$\pi_*: R^*(\overline{M}_{g_1,n_1\scup\{*\}}) 
\otimes_{\Q}
R^*(\overline{M}_{g_2,n_2\scup\{\bullet\}}) \rarr
R^*(\overline{M}_{g_1+g_2, n_1+n_2}),$$
$$\pi_*: R^*(\overline{M}_{g_1, n_1\scup\{*,\bullet\}}) \rarr
R^*(\overline{M}_{g_1+1, n_1}).$$
\end{enumerate}

The tautological rings  possess
remarkable algebraic 
and combinatorial structures with basic connections
to topological gravity. A discussion of these properties together
with a conjectural
framework for the study of $R^*(\overline{M}_{g,n})$
can be found in [F], [FP].

In genus 0, the equality
\begin{equation*}
R^*(\overline{M}_{0,n}) = A^*(\overline{M}_{0,n}),
\end{equation*}
for  $n\geq 3$, is well-known from Keel's study [K].

Denote the image of
$R^*(\overline{M}_{g,n})$ under the canonical 
map to the ring of {\em even} cohomology
classes by:
$$RH^*(\overline{M}_{g,n})\subset H^{2*}(\overline{M}_{g,n}).$$
In genus 1, Getzler has claimed the isomorphisms:
$$ R^*(\overline{M}_{1,n}) \eqq RH^*(\overline{M}_{1,n}),$$
\begin{equation*}
RH^*(\overline{M}_{1,n}) = H^{2*}(\overline{M}_{1,n}),
\end{equation*}
for $n\geq 1$, see [G1].

For $g>1$, complete results are known only in codimension 1. The equality
$$R^1(\overline{M}_{g,n}) = A^*(\overline{M}_{g,n}),$$
for $2g-2+n>0$, is a consequence of Harer's cohomological
calculations [H]. 

It is natural to ask whether all algebraic cycles classes
on $\M_{g,n}$ are tautological. The existence of nontautological cycles
defined over $\com$ may be deduced from the odd cohomology of
$\overline{M}_{1,11}$. There are two arguments which may be used:
\begin{enumerate}
\item[(i)]
By a theorem of Jannsen, since the map to cohomology
$$A^*(\overline{M}_{1,11}) \rarr H^*(\overline{M}_{1,11})$$
is not surjective, the map is not injective. We may then deduce the
existence of a nontautological Chow class in $A^*(\overline{M}_{1,11})$
from Getzler's claims \cite{b}.

\item[(ii)] More precisely, 
the existence of a holomorphic 11-form and
a theorem of Srinivas together imply  
$A_0(\M_{1,11})$ 
is an infinite dimensional vector space,  
 while 
$R_0(\M_{1,11}) \eqq \Q$, see \cite{grv}. 
\end{enumerate}

These arguments do not produce an explicit algebraic cycle which is
not tautological. Several further questions are also left open.
Are there nontautological cycles defined over 
$\Q$?  Are there 
algebraic cycles with cohomological image not
contained in $RH^*(\overline{M}_{g,n})$?  Are there
nontautological classes on the noncompact spaces $M_{g,n}$?

We answer all these questions in the affirmative by explicit
constructions of 
integrally defined algebraic
cycles.  Our basic
criterion for detecting nontautological cycles is the following 
Proposition.

\begin{pr}\label{criterion}
Let $\iota : \overline{M}_{g_1,n_1\cup \{*\}} \times 
\overline{M}_{g_2,n_2\cup\{\bullet\}} 
\rarr \overline{M}_{g,n_1+n_2}$
be the gluing map to a boundary divisor.
If $\gamma \in RH^*(\overline{M}_{g_1+g_2,n_1+n_2})$, 
then $\iota^*(\gamma)$ has
a tautological K\"unneth decomposition:
$$\iota^*(\gamma) \in 
RH^*(\overline{M}_{g_1,n_1\cup\{*\}}) \otimes 
RH^*(\overline{M}_{g_2,n_2\cup\{\bullet\}}).$$
\end{pr}
While the above result
is well known to experts, we know of no adequate reference,
so we give a proof in the Appendix.

Our strategy for finding nontautological classes 
combines Proposition \ref{criterion} with the
existence of odd cohomology on the moduli spaces of curves.
We find loci in moduli
space which restrict to diagonal loci of symmetric
boundary divisors.
By the existence of odd cohomology in certain cases,
the K\"unneth decomposition of the diagonal is not tautological.

Let $h$ be an {\em odd} integer and set $g=2h$.  Let $Y\subset \M_g$ 
denote the closure of the set of nonsingular curves of genus $g$ 
which admit a degree 2 map to a nonsingular curve of genus $h$. 
Intersecting $Y$ with the boundary map from 
$\M_{h,1} \times \M_{h,1}$ yields the diagonal. 
Pikaart has proven, for sufficiently large $h$,
$\M_{h,1}$ has odd cohomology [P].
Hence, we can conclude $Y$ is not a tautological class,
even in homology.

\begin{tm}\label{gbar} For all sufficiently large odd $h$,
$$[Y] \notin RH^*(\overline{M}_{2h}).$$
\end{tm}

Our other examples are loci in the moduli space of
pointed genus 2 curves. We will use the
 odd cohomology of $\M_{1,11}$ to find nontautological
K\"unneth decompositions.

Let $\sigma$ in $\bbS_{20}$ be a product of 10 disjoint 2-cycles,
$$\sigma= (1, 11)(2, 12) \cdots (10, 20),$$
inducing an involution on $\M_{2,20}$.  Let
$Z$ denote the component of the fixed locus of the involution
corresponding generically to a 20-pointed, nonsingular,
 bielliptic curve of genus 2 with
the 10 pairs of conjugate markings.  
$Z$ is of codimension 11 in $\M_{2,20}$.
The intersection of $Z$ 
with the boundary map
\begin{equation}
\label{kqw}
\iota:\M_{1,11} \times \M_{1,11} \rarr \M_{2,20},
\end{equation} 
yields the diagonal.  

\begin{tm}\label{2bar} 
$[Z] \notin RH^*(\overline{M}_{2,20}).$
\end{tm}

Although the methods used to prove Theorems \ref{gbar} and \ref{2bar}
depend crucially on the structure
of the boundary of the moduli space, in Section \ref{2nobarsec}
we use Getzler's results on the cohomology of $\M_{1,n}$
to show the class $[Z]$ is nontautological even on the interior.

\begin{tm}\label{2nobar}  $[Z] \notin R^*(M_{2,20}).$
\end{tm}

Finally, although the diagonal loci were used in our deductions of
the above results, we could not conclude the diagonals themselves
were nontautological.  We show a diagonal locus is
nontautological in at least one case.
Let $\iota$ denote the boundary inclusion,
$$\iota: \overline{M}_{1,12} \times \overline{M}_{1,12} \rarr
\overline{M}_{2,22},$$
Let $\Delta$ denote the
class of the diagonal in 
$A^*(\overline{M}_{1,12} \times \overline{M}_{1,12})$.

\begin{tm}\label{diagonal} The push-forward $\iota_*[\Delta]$ is not a 
tautological class:
$$\iota_*[\Delta] \notin RH^*(\overline{M}_{2,22}).$$
\end{tm}

\noindent While it seems likely the image of the diagonal 
by (\ref{kqw})
in $\M_{2,20}$ is not tautological, we do not have a proof.

To our knowledge
there are still no (proven) examples of nontautological classes on $M_g$.  
While the methods 
of our paper could perhaps be used to find such a class, in particular
the class of Theorem 1 may be nontautological when
restricted to $M_g$, 
our techniques are unlikely to produce nontautological 
classes of low codimension.
As the tautological ring of $M_g$ vanishes in codimension $g-1$ and higher,
the question of  nontautological classes
on $M_g$ of codimension less than $g-1$ is particularly interesting.

The authors thank C. Faber and E. Getzler for several 
conversations about the cohomology
of $\overline{M}_{1,n}$, and R. Vakil for discussions about
hunting for nontautological classes.
 T.~G.~ was partially supported by an NSF
post-doctoral fellowship. R.~P.~ was partially supported by DMS-0071473 and
fellowships from the Sloan and Packard foundations. Part of the
research reported here 
was pursued during a visit by R.~P.~ to the Royal Technical
Institute in Stockholm.

\section{Admissible double coverings}

For the proofs of Theorems \ref{gbar} and \ref{2bar} we will require
certain moduli spaces of double covers.  Choose
$g$ and $h$ with $g\geq 2h-1$.  We let $M(g,h)$ denote
the (open) space parameterizing {\em double} covers,
$$\pi:C_g \rarr C_h,$$
of curves of genus $g$ and $h$ respectively 
{\em together} with an ordering of
the branch points of the morphism $\pi$.  
The space $M(g,h)$ is a finite \'etale cover of
$M_{h,b}$ where $b=2(g-2h+1)$ is the number of branch points of
$\pi$.  
The map,
$$\mu:M(g,h) \rarr M_{h,b},$$ is simply 
defined by
$$\mu([\pi]) = [C_h, p_1, \ldots, p_b],$$
where $p_1, \ldots, p_b$ are the ordered 
branch points.  

There is a natural compactification by admissible double covers, 
$$M(g,h) \subset \M(g,h),$$ over $\M_{h,b}$.  An admissible
double cover $\pi$ of a stable curve
is branched over the marked points and possibly the nodes.
Over the nodes of the target, 
the map $\pi$ is either \'etale or \'etale locally
of the form,
$$\pi: Spec(\com[x,y]/(xy)) \rarr Spec(\com[u,v]/(uv)),$$
$$ u=x^2, v=y^2.$$
By construction,
the space $\M(g,h)$ 
is equipped with maps to both $\M_{h,b}$ and $\M_g$.  
The latter map involves a stabilization process, since the source
curve of an admissible covering need not be stable.

We will require additionally pointed
 moduli spaces of admissible covers,
$\M_k(g,h)$.  These pointed spaces are  finite covers of
$\M_{h,b+k}$ which parametrizes admissible double covers of a 
$b+k$ pointed
nodal curve of genus $h$ by a curve of genus $g$, with the ramification
over the first $b$ marked points and possibly 
the nodes of the target curve together
with an ordering of the fibers of the last $k$ marked points.  The
pointed
spaces 
are equipped with natural morphisms to $\M_{h,b+k}$ and $\M_{g,2k}$.
For the latter map, we adopt the ordering convention that the two
points in the fiber over the $(b+i)$th marked point of the target curve
have markings $i$ and $k+i$ on the source.

Essentially, we require only one fact about the moduli spaces
of admissible covers:
$M(g,h)\subset \M(g,h)$ is dense (and similarly for the open subset
$M_k(g,h) \subset \M_k(g,h)$).  
Over the complex numbers, the density 
is easily proven analytically.  One can simply locally smooth
the double cover of a small neighborhood of the node, and glue the 
result together with the restriction of the original cover away 
from the node.  The local description of our double covers ensures
that they can be smoothed locally.  A treatment of the
theory of admissible covers can be found in [HaM].

\section{Proof of Theorems \ref{gbar} and \ref{2bar}}

Let $h$ be an odd positive integer, and let $g=2h$.  
Consider the morphism 
$$\phi:\M(g,h) \to \M_g.$$
The image cycle, 
$$Y=\phi(\M(g,h)),$$ consists of those curves of genus
$g$ which are admissible double covers of a curve of genus $h$.  
Equivalently, $Y$ is the closure of the set of nonsingular curves
of genus $g$ which admit a degree 2 map to a nonsingular curve
of genus $h$.  We want
to apply Proposition \ref{criterion} to conclude that $[Y]$ is
not tautological.  We will look at the pullback of $[Y]$ under
the gluing map
$$\iota:\M_{h,1} \times \M_{h,1} \to \M_g.$$ 

\begin{lm} \label{joq}
$\iota^*([Y]) = c [\Delta ]$ for some positive constant $c$.
\end{lm}

\bpf

We first prove $\Delta \subset \iota^{-1}(Y)$. Let $[C,p]\in 
\overline{M}_{h,1}$.
We will construct an 
admissible double cover with target $C$ union a
a rational tail glued at $p$ carrying the two branch markings.  
A double cover is given by two disjoint 
copies of $C$ joined by a rational curve with a degree 2 mapping  to the
rational tail of the target branched over the two markings.  Under
stabilization, the domain is mapped to the diagonal point
$$\iota( [C,p] \times [C,p]).$$
An easy count shows $\Delta$ is an irreducible 
component
of $\iota^{-1}(Y)$ of expected dimension.  

To prove the Lemma, we need only show
$\Delta = \iota^{-1}(Y)$.
Suppose there were another irreducible component $I$.  Let
$$\pi:C_g \rarr C_h$$
 be an admissible double cover corresponding to a general 
point of $I$.  Then $C_g$ may be expressed
as a union of two curves of arithmetic
genus $h$ joined at a single node. 
The chosen node of $C_g$ must map to a node of $C_h$.
Since the space of admissible
coverings is a finite cover of $\M_{h,2}$, the preimage of the locus of 
curves with 2 or more nodes is not a divisor.  Hence, we conclude 
$C_h$ has {\em exactly} 1 node. 

 The  node of $C_h$ must be disconnecting since there
are no reducible admissible double covers of an irreducible curve with 
branch points.  We write $C_h=T_1 \cup T_2$.  Since $h$ is odd,
we may assume, without loss of generality,
$T_1$ has genus greater than $h/2$. 

Since $C_h$ has 1 node, $C_g$ must have either
1 or 2 nodes. Since
any cover of $T_1$ 
by a curve of genus $h$ must be unramified,
$C_h$ cannot have exactly 1 node.  

The domain $C_h$ must therefore have 
2 nodes lying over the node of $C_g$.
If the induced
cover of $T_1$ were connected, the neither node of $C_g$
could be disconnecting  
Hence, the cover of $T_1$ must be disconnected.

Then, each component of the cover of $T_1$ must map
isomorphically to $T_1$. The cover of $T_2$ must be
connected of genus 0 in order for the assumed
decomposition of $C_g$ into curves of arithmetic genus $h$
to exist. Therefore, we find we are in the component
$\Delta$ of $\iota^{-1}(Y)$.
\epf

Pikaart \cite{pi} has shown 
for all sufficiently large values of $h$, $$H^{33}(\M_{h,1}) \neq 0.$$  
Hence, the diagonal in $\M_{h,1} \times \M_{h,1}$
does not have tautological K\"unneth decomposition.
By Proposition \ref{criterion}, the proof of
Theorem \ref{gbar} is complete.

The argument for the nontautological cycle on $\M_{2,20}$ is
similar.  Let $Z$ be the
image of $\M_{10}(2,1)$ in $\M_{2,20}$.
Consider the boundary stratum,
$$\iota:\M_{1,11} \times \M_{1,11} \to \M_{2,20},$$
obtained by attaching at the 
last point on each marked curve, and
numbering the markings of the glued curve in order with the first
10 markings from the first factor and the last 10 from the
second factor.

\begin{lm} $\iota^*([Z]) = c [\Delta]$ for some positive constant $c$.
\end{lm}
The proof of the Lemma is essentially identical to the proof of
Lemma \ref{joq} above.  Theorem \ref{2bar} is then
a consequence of Proposition \ref{criterion} and 
the existence of odd cohomology on $\M_{1,11}$.

\section{Proof of Theorem \ref{2nobar}} \label{2nobarsec}

To deduce Theorem \ref{2nobar} from Theorem \ref{2bar}, 
we will need the following  
results announced by Getzler:
\begin{equation}
\label{kqp}
RH^*(\overline{M}_{1,n}) = H^{2*}(\overline{M}_{1,n}),
\end{equation}
and for all odd $k<11$,
\begin{equation}
\label{kpp}
H^k(\overline{M}_{1,n})=0.
\end{equation}
The statement (\ref{kqp}) 
is equivalent to the generation of even cohomology by
the classes of boundary strata for $\overline{M}_{1,n}$.
Actually, we require the following consequences of
Getzler's results.

\begin{lm}\label{ezracor}
Three properties for $\overline{M}_{1,n}$:

\begin{enumerate}
\item[(i)] Every algebraic cycle on 
$\M_{1,11} \times \M_{1,11}$ 
of complex codimension less than 11 is
homologous to a tautological class.
\item[(ii)] 
Every algebraic cycle on 
$\M_{1,m} \times \M_{1,n} \times \prod_i \M_{0,{l_i}}$ 
is homologous to a tautological class for $m<11$.
\item[(iii)] Every algebraic cycle on $\M_{1,n} \times \prod_i \M_{0,l_i}$
is homologous to a tautological class.
\end{enumerate}
\end{lm}

\bpf
Let
$V$ be an algebraic cycle on $\M_{1,11} \times \M_{1,11}$ 
of complex codimension less than 11. 
 Consider the K\"unneth decomposition
of $[V]$.  There can be no odd terms by (\ref{kpp}).
Thus, by (\ref{kqp}), we can write
$[V]$ as a sum of products of tautological classes proving (i).
By (\ref{kpp}) and
 Poincar\'e duality, all the cohomology of $\M_{1,m}$ is
tautological when $m<11$.
By Keel's results, 
all the cohomology of $\M_{0,l_i}$ is tautological.  
Hence,
in the K\"unneth decomposition of our cycle in parts (ii) and (iii), 
none of the odd
cohomology of $\M_{1,n}$ can appear. 
\epf

Consider the class $[Z]$ on $\M_{2,20}$ constructed in Theorem 
\ref{2bar}.
We claim the image of $[Z]$ in $A^*(M_{2,20})$ is not tautological.  
The argument is by contradiction.

Suppose the image
is tautological. 
There must exist a collection of cycles
$Z_i$ of codimension equal to  11 in $\M_{2,20}$, 
supported on boundary strata, for which $Z + \sum Z_i$
is tautological.  Then, $\sum Z_i$ is not
homologous to a tautological class when intersected with 
$\M_{1,11} \times \M_{1,11}$. 

By
Lemma \ref{ezracor} part (i), if any cycle
 $Z_i$ is supported on the image stratum of  $\M_{1,11} \times
\M_{1,11}$, then $Z_i$ is homologous to a
tautological class (since the codimension of $Z_i$ is
less than 11 in the divisor).
We discard all $Z_i$ contained in the image of
$\M_{1,11} \times
\M_{1,11}$.

Let $X$ be the union of boundary divisors supporting
the remaining $Z_i$. The sum of the remaining $Z_i$
is homologically nontautological when pushed into 
$\M_{2,20}$ and restricted to $\M_{1,11} \times \M_{1,11}$.
However,
it is clear the push-pull will produce an algebraic cycle class
supported on 
\begin{equation}
\label{njk}
X\ \ \cap\ \  \M_{1,11} \times \M_{1,11}.
\end{equation}
Since $X$ does not contain the image of $\M_{1,11} \times \M_{1,11}$,  
the locus  (\ref{njk}) is contained
in  boundary strata which either 
have a genus 1 factor with fewer than 11 points,
or have less than two genus 1 factors.
Parts (ii) and (iii) of Lemma \ref{ezracor} show that there are no
homologically nontautological classes supported on these loci. 
The contradiction completes  the proof
of Theorem \ref{2nobar}.

\section{Proof of Theorem \ref{diagonal}}

\subsection{Odd cohomology of $\M_{1,n}$}

We will require several properties of  the odd cohomology of the
moduli spaces $\overline{M}_{1,n}$ for the proof of Theorem \ref{diagonal}. 
The first is a well-known 
specialization of (\ref{kpp}).
\begin{pr} \label{vanishing}
The odd cohomology groups of $\overline{M}_{1,n}$ vanish
 in case $1\leq n < 10$.
\end{pr}

Cusp forms of weight $n$ may be used to construct cohomology classes
in $H^{n-1,0}(\overline{M}_{1,n-1}, \com)$. The discriminant
form $\Delta$, the unique cusp form of weight 12, yields a canonical
non-zero element $s \in H^{11,0}(\overline{M}_{1,11}, \com)$. 

\begin{pr}\label{odd11} 
The odd cohomology of $\overline{M}_{1,11}$ is concentrated
          in $$H^{11,0}(\overline{M}_{1,11},\com)\eqq\com,$$ 
$$H^{0,11}(\overline{M}_{1,11},\com)\eqq\com.$$
          Moreover, the ${\mathbb{S}}_{11}$-module
          in both cases is the alternating representation.
\end{pr}
\noindent By the ${\mathbb{S}}_{11}$-module identification, the class
$s$ is {\em not} ${\mathbb{S}}_{11}$-invariant.
Let $t\in H^{0,11}(\overline{M}_{1,11}, \com)$
denote the uniquely defined Poincar\'e dual class to $s$:
$$\int_{\overline{M}_{1,11}} s \cup t = 1.$$

Both Proposition \ref{vanishing} and \ref{odd11} are well-known. 
Proofs can be found,
for example, in [G2] where the ${\mathbb{S}}_n$-equivariant Hodge polynomials
of $\overline{M}_{1,n}$ are calculated for all $n$. 
We will need a dimension calculation in the  $n=12$ pointed case [G2]:
\begin{pr}\label{dim12} 
The dimension of $H^{11,0}(\overline{M}_{1,12}, \com)$ is 11.
\end{pr}
\noindent In fact, the odd cohomology of $\overline{M}_{1,12}$ is
concentrated in $H^{11,0}$, $H^{0,11}$, $H^{12,1}$, and $H^{1,12}$
(all 11 dimensional).

\subsection{ A basis for $H^{11,0}(\overline{M}_{1,12}, \com)$}
Let $S=\{1,2,3,\ldots, 11, p\}$. 
For each index $1\leq i \leq 11$, let
$$\pi_i: \overline{M}_{1,S} \rarr \overline{M}_{1, S- i} \eqq 
\overline{M}_{1,11}$$
denote the forgetful map.
Since we consider $$S-i = \{1,2,3, \ldots, \hat{i}, \ldots,11, p\}$$ as
an {\em ordered} set, the last isomorphism above is canonical.
Define the classes $a_i$ and $b_i$ by:
 $$a_i = \pi_i^*(s) \in H^{11,0}(\overline{M}_{1,S}, \com),$$
 $$b_i = \pi_i^*(t) \in H^{0,11}(\overline{M}_{1,S}, \com).$$

For each index $1 \leq i \leq 11$, let $\epsilon_i$ be the map
defined by the inclusion
$$\epsilon_i: \overline{M}_{1,11} \eqq 
\overline{M}_{1, S- i} \rarr \overline{M}_{1,S}.$$
Here, an $S$ pointed curve is obtained from an
$S-i$ pointed curve by attaching a rational tail
containing the markings $i,p$ to the point $p$ of the latter curve.
The map $\epsilon_i$ is simply the inclusion of the boundary divisor
$D_{ip}$ with genus splitting $1+0$ and point splitting
$$\{1, \ldots, \hat{i}, \ldots, 11\} \ \cup\  \{i,p\}.$$
Define the classes $c_i$ and $d_i$ by:
 $$c_i = \epsilon_{i*}(s) \in H^{12,1}(\overline{M}_{1,S}, \com),$$
 $$d_i = \epsilon_{i*}(t) \in H^{1,12}(\overline{M}_{1,S}, \com).$$
Here, the  cohomological push-forward is defined by the
equivalent equalities
$$\int_{\overline{M}_{1,S}} \epsilon_{i*}(x) \cup y =
\int_{\overline{M}_{1,11}} x \cup \epsilon_i^*(y),$$
$$
\int_{\overline{M}_{1,S}} y \cup \epsilon_{i*}(x)  =
\int_{\overline{M}_{1,11}} \epsilon_i^*(y) \cup x.$$

\begin{pr}\label{intform1} The sets $\{a_1, \ldots, a_{11}\}$ and
$\{d_1, \ldots, d_{11}\}$ form a pair of Poincar\'e dual bases of
$H^{11,0}(\overline{M}_{1,S}, \com)$ and $H^{1,12}(\overline{M}_{1,S}, \com)$.
\end{pr}

\bpf
By the dimension result of Proposition \ref{dim12}, it suffices to prove
\begin{equation}
\label{nnsd}
\int_{\overline{M}_{1,S}} a_i \cup d_j = \delta_{ij}.
\end{equation}
By definition of the cohomological push-forward,
\begin{equation}
\label{sdf}
\int_{\overline{M}_{1,S}} \pi^*_i(s) \cup \epsilon_{i*}(t) =
\int_{\overline{M}_{1,11}} s\cup t =1.
\end{equation}
The first equality in (\ref{sdf}) is true exactly (not up to sign)
by the precise ordering conventions used.

The vanishing of (\ref{nnsd}) when $i\neq j$ is a direct
consequence of Proposition \ref{vanishing}. We find
\begin{equation}
\label{sdff}
\int_{\overline{M}_{1,S}} \pi^*_i(s) \cup \epsilon_{j*}(t) =
\int_{\overline{M}_{1,11}} \epsilon_j^* \pi^*_i(s)\cup t.
\end{equation}
The composition $\pi_i \circ \epsilon_j $ has image
isomorphic to $\overline{M}_{1,10}$. As the image 
supports no odd cohomology, the integral (\ref{sdff})
vanishes.
\epf

An identical argument proves the duality result for the
classes $c_i$ and $b_i$.

\begin{pr}\label{intform2} The sets $\{c_1, \ldots, c_{11}\}$ and
$\{b_1, \ldots, b_{11}\}$ form a pair of Poincar\'e dual bases of
$H^{12,1}(\overline{M}_{1,S}, \com)$ and $H^{0,11}(\overline{M}_{1,S}, \com)$.
The intersection form is:
\begin{equation*}
\int_{\overline{M}_{1,S}} c_i \cup b_j = \delta_{ij}.
\end{equation*}
\end{pr}


\subsection{The action of $\psi_p$}
Let $\psi_p\in H^{1,1}(\overline{M}_{1,S}, \com)$ denote the cotangent line
class at the point $p$. Multiplication by $\psi_p$ defines
linear maps:
$$ \Psi: H^{11,0}(\overline{M}_{1,S}, {\com}) \rarr H^{12,1}(\overline{M}_{1,S}, \com),$$
$$ \Psi: H^{0,11}(\overline{M}_{1,S}, {\com}) \rarr H^{1,12}(\overline{M}_{1,S}, \com).$$
These maps are completely determined by the following result.

\begin{pr}\label{psiaction} For all $1\leq i \leq 11$,
$\Psi(a_i) = c_i$ and $\Psi(b_i) = d_i$.
\end {pr}

\bpf Consider the morphism $\pi_i: \overline{M}_{1,S} \rarr \overline{M}_{1,S-i}$.
A standard comparison result governing the cotangent line class is:
$$\psi_p = \pi_i^*(\psi_p) + [D_{ip}],$$
where $\pi_i^*(\psi_p)$ denotes the pull-back of the cotangent class
on $\overline{M}_{1, S-i}$.
We then find:
\begin{equation}
\label{llssl}
\psi_p \cup \pi_i^*(s) = \pi_i^*(\psi_p \cup s) + [D_{ip}] \cup \pi_i^*(s).
\end{equation}
Since $\overline{M}_{1,11}$ has odd cohomology only in degree 11, the
first summand of (\ref{llssl}) vanishes. The second summand is exactly
equal to $\epsilon_{i*} (s)$ (using the ordering conventions).
We conclude $\Psi(a_i)=c_i$.
The derivation of $\Psi(b_i)=d_i$ is identical.
\epf

\subsection{Proof of Theorem \ref{diagonal}}

Set $\tilde S=\{\tilde 1, \tilde 2, \ldots, \tilde 11, \tilde p\}$.  
Consider the boundary map 
$$\iota: \overline{M}_{1,S} \times
\overline{M}_{1,\tilde{S}} \rarr \overline{M}_{2,22}$$ 
defined by attaching $p$ to $\tilde p$ (and ordering the markings
arbitrarily). Define $\gamma$ by
$$\gamma = \iota_*[\Delta] \in H^{26}(\overline{M}_{2,22}),$$
where $\Delta$ is the diagonal subvariety of 
$\overline{M}_{1,S} \times
\overline{M}_{1,\tilde{S}}$ (under the canonical
isomorphism $\overline{M}_{1,S} \eqq \overline{M}_{1,\tilde{S}}$).

Here, $\iota$ is easily seen to define an {\em embedding}.
The normal bundle to $\iota$ in $\overline{M}_{2,22}$ has top Chern class
$-\psi_p - \psi_{\tilde{p}}$. By the self-intersection formula,
$$\iota^* \iota_*[\Delta] = [\Delta] \cup (-\psi_p -\psi_{\tilde{p}}).$$

Let 
$X_1, \ldots, X_m$ be a basis of $H^*(\overline{M}_{1,S})$.
Let $\tilde{X}_1, \ldots, \tilde{X}_m$ denote the corresponding
basis of $H^*(\overline{M}_{1,\tilde{S}})$. The K\"unneth decomposition
of $[\Delta]$ 
 is determined by:
$$[\Delta] = \sum_{i,j} g^{ij} X_i \otimes \tilde{X}_j \in
H^*(\overline{M}_{1,S}) \times \overline{M}_{1,\tilde{S}}),$$
where 
$$g_{ij} = \int_{\overline{M}_{1,S}} X_i \cup X_j.$$
In particular, if $X_1, \ldots, X_m$ is a self dual basis, then
$$[\Delta] = \sum_{i} (-1)^{\nu_i \nu^\vee_i} X_i \otimes X^\vee_i,$$
where $\nu_i, \nu_i^\vee$ are the degrees of $X_i$ and $X_i^\vee$
respectively.

We are interested in the K\"unneth components of 
$[\Delta]$ of odd type --- that is K\"unneth components
lying in  $$H^{odd}(\overline{M}_{1,S}) \otimes 
H^{odd}(\overline{M}_{1,\tilde{S}}).$$
By Proposition \ref{intform1} and \ref{intform2},
the odd type summands of $[\Delta]$ are:
$$ \sum_{i=i}^{11} - a_i \otimes \tilde{d_i} +b_i \otimes \tilde{c_i}
-c_i \otimes \tilde{b_i} + d_i \otimes \tilde{a_i}.$$
Hence, the odd summands of $\iota^*\iota_*[\Delta]$ are:
$$ \sum_{i=1}^{11} \Psi(a_i) \otimes \tilde{d_i} - \Psi(b_i) \otimes \tilde{c_i}
+ c_i \otimes \tilde{\Psi} (\tilde{b_i}) - d_i \otimes \tilde{\Psi}(\tilde{a_i}).$$
By Proposition \ref{psiaction}, we find the odd summands of 
$\iota^*\iota_*[\Delta]$
equal:
\begin{equation}
\label{finny}
\sum_{i=1}^{11} 2 c_i \otimes \tilde{d_i} - 2 d_i \otimes \tilde{c_i}.
\end{equation}
As the odd summands (\ref {finny}) do not vanish,
$$\iota_*[\Delta] \notin RH^*(\overline{M}_{2,22})$$
by Proposition \ref{criterion}.
The proof of Theorem \ref{diagonal} is complete.
\epf

\appendix
\section{Pull-backs in the tautological ring}
\subsection{Stable graphs}
The boundary
strata of the moduli space of curves correspond
to {\em stable graphs} $$A=(V, H,L, g:V \rarr {\mathbb Z}_{\geq 0}, a:H\rarr V, i: H\rarr H)$$
satisfying the following properties:
\begin{enumerate}
\item[(i)] $V$ is a vertex set with a genus function $g$,
\item[(ii)] $H$ is an half-edge set equipped with a 
vertex assignment $a$ and fixed point free involution
$i$,
\item[(iii)] $E$, the edge set, is defined by the
orbits of $i$ in $H$ (self-edges at vertices
are permitted),
\item[(iv)] $(V,E)$ define a {\em connected} graph,
\item[(v)] $L$ is a set of numbered legs attached to the vertices,
\item[(vi)] For each vertex $v$, the stability condition holds:
$$2g(v)-2+ n(v) >0,$$
where $n(v)$ is the valence of $A$ at $v$ including both half-edges and legs.

\end{enumerate}
The genus of $A$ is defined by:
$$g(A)= \sum_{v\in V} g(v) + h^1(A).$$
Let $v(A)$, $e(A)$, and $n(A)$ denote the cardinalities of $V$, $E$, and $L$
respectively.
A boundary stratum of $\overline{M}_{g,n}$ naturally determines
a stable graph of genus $g$ with $n$ legs by considering the dual graph of 
a generic pointed curve parameterized by the stratum.

Let $A$ be a stable graph. Define the moduli space
$\M_A$ by the product: 
$$\M_A =\prod_{v\in V(A)} \M_{g(v),n(v)}.$$
 Let $\pi_v$ denote the projection from $\M_A$ to 
$\M_{g(v),n(v)}$ associated to the vertex $v$.  There is a
canonical
morphism $\xi_{A}: \M_{A} \rarr \M_{g,n}$ with image equal to the boundary stratum
associated to the graph $A$.  To construct $\xi_A$, 
a family of stable pointed curves over $\M_A$ is required.  Such a family
is easily defined 
by attaching the pull-backs of the universal families over each of the 
$\M_{g(v),n(v)}$  along the sections corresponding to half-edges.

\subsection{Specialization}
\label{specl}
Our main goal in the Appendix is to understand the fiber product:
$$
\begin{CD} F_{A,B} @>>> \M_B \\
@VVV @V{\xi_B}VV\\
\M_A @>{\xi_A}>> \M_{g,n}. \end{CD}
$$
To this end, we will require additional terminology.  
A stable graph $C$ is a {\em specialization} of a stable graph $A$
if $C$ is obtained from $A$ by replacing each 
vertex $v$ of $A$ with a stable graph of genus $g(v)$ with $n(v)$ legs. 
Specialization of graphs corresponds to
specialization of stable curves.  

There is a subtlety involved in the notion of specialization: 
a given graph $C$ may arise as a specialization of $A$ in more than one
way.  An {\em $A$-graph structure} on a stable graph
$C$ is a choice of subgraphs of $C$ in bijective
correspondence with $V(A)$ such that
$C$ can
be constructed by replacing each vertex of $A$ by the corresponding
subgraph.  If $C$ has an $A$-structure, 
then every half-edge of $A$ corresponds to a particular
half-edge of $C$, and every vertex of $C$ is associated to
a particular vertex of $A$.  

A point of $\M_A$ is given by a stable curve together with a choice of
$A$-structure on its dual graph.  In fact, we can naturally
identify the stack $\M_A$ with a stack defined in terms
of $A$-structures. This identification will be useful for the
analyzing the fiber products of strata.

Define a stable $A$-curve over a connected base $S$,
$$\pi:\CC \rarr S,$$ to be
a stable $n(A)$-pointed curve of genus 
$g(A)$ over $S$ together with:
\begin{enumerate}
\item[(i)] $e(A)$ sections $\sigma_1, \ldots, \sigma_{e(A)}$
of $\pi$ with image in the singular
locus of $\CC$, 
\item[(ii)] $2e(A)$ sections of the normalization of
$\CC$ along the sections $\{\sigma_i\}$ corresponding to the nodal separations,
\item[(iii)] $v(A)$ disjoint 
$\pi$-relative components of $\CC \setminus \{\sigma_i\}$ whose union
is $\CC \setminus \{ \sigma_i \}$,
\item[(iv)]
An isomorphism between $A$ and the canonical
stable graph defined by the dual graph of the $v(A)$ $\pi$-relative
components and $2e(A)$ sections of the normalization (corresponding to half-edges).
\end{enumerate}
Here,  a $\pi$-relative component is a connected component of 
of $\CC\setminus \{\sigma_i\}$ which remains connected upon
pullback under an arbitrary morphism of connected schemes $h:T\rarr S$.

The data of a stable $A$-curve can be pulled back 
under any morphism of base schemes.
After pull-back to a geometric point, an $A$-curve is exactly 
an $A$-structure on the dual graph of the corresponding curve.  

A stack $\M_A'$ of curves with $A$-structure 
morphisms and respecting the $A$-structure may be defined. However, we 
find the following result.

\begin{pr}There is a natural isomorphism between $\M_A$ and $\M_A'$.
\end{pr}  
\bpf
A natural morphism from $\M_A$ to $\M_A'$ is obtained by  
assigning the canonical $A$-structure to the universal curve over
$\M_A$.
In the other direction, given
an $S$-valued point of $\M_A'$, we naturally obtain a collection of $v(A)$ 
stable curves by analyzing the $\pi$-relative components of $C$ normalized
at the $e(A)$ nodes.  Since we have a bijection
between these curves and $v(A)$, and a bijection between the new markings
and the $2e(A)$ sections, we obtain an $S$-valued point of $\M_A$.  
This correspondence induces a bijection on the space of morphisms
between corresponding objects.
\epf

\subsection{Fiber products}
By definition,
an $S$-valued point of $F_{A,B}$
is an $S$-valued point of $\M_A$, an $S$-valued point of $\M_B$, and a choice
of isomorphism between the two 
pull-backs of the universal curve over $\M_{g,n}$
under the boundary inclusions.  If $S$ is $\Spec (\com)$,
we find  the dual graph $C$ of the curve over $S$
defined by the map to $\M_{g,n}$ is naturally equipped with both an 
$A$-structure and a $B$-structure.  Conversely, given a curve $C$
together with two such structures on the dual graph, we 
naturally obtain a point
of $F_{A,B}$.  A graph $C$ equipped with both $A$ and $B$-structures will
be called 
an $(A,B)$-graph.

An $(A,B)$-graph $C$  
is  {\em generic} if every half-edge of $C$ corresponds to a half-edge of $A$
or a half-edge of $B$.
The irreducible components of $F_{A,B}$ will
correspond to generic $(A,B)$-graphs.
A graph with an $(A,B)$-structure is canonically
a specialization of a unique generic $(A,B)$-graph: the generic graph
is obtained by contracting all those edges which do not correspond to 
edges of $A$ or $B$.

Associated to an $(A,B)$-graph $C$, we obtain a moduli space
$\M_C$ which naturally maps to $F_{A,B}$.  The moduli space may be described
either as $\prod_{v\in v(C)} \M_{g(v),n(v)}$ or in stack terms analogous to the
definition of $\M_A'$ above (the stack does not depend on the 
$(A,B)$-structure on $C$, although the map to $F_{A,B}$ does).  
We find the following result.

\begin{pr} \label{freddy}
There is a canonical isomorphism between $F_{A,B}$ and the disjoint
union of $\M_C$ over all generic $(A,B)$-graphs $C$.
\end{pr}

\bpf
It will suffice to identify the categories involved over
connected base schemes $S$.  We will give the morphisms in both directions.

If $C$ is an $(A,B)$-graph, then we clearly have a morphism from $\M_C$
to both $\M_A$ and $\M_B$, and a choice of isomorphism between the induced
maps to $\M_{g,n}$.  

In the other direction, suppose data corresponding to $F_{A,B}$ is given
over $S$.
In particular, we have a stable curve over $S$:
$$\pi: {\mathcal C} \rarr S.$$  Consider the $\pi$-fiber
over a geometric point of $S$.  The $\pi$-fiber has a dual graph which
is equipped with an $(A,B)$-structure by virtue of the maps
to $\M_A$ and $\M_B$.  Let $C$ be the unique generic $(A,B)$-graph which
specializes to the $(A,B)$-structure found at the geometric point. 

There is a canonical $C$-structure on $\mathcal{C}$. 
The half-edges of $C$
are already naturally identified with half-edges of the graphs $A$ or $B$, and
since ${\mathcal {C}}$ has $A$ and $B$-structures, $\pi$ is  equipped with
sections associated to all of the half-edges. 
The $C$-structure is constant on ${\mathcal C}$ because of the
connectedness of $S$.

The morphisms in the two categories are the same by a 
straightforward check. However,  it is important to note that an automorphism
of an object of $F_{A,B}$ must induce a trivial automorphism of the
graph $C$, because each half-edge corresponds to an edge of either $A$ or $B$.
\epf

\subsection{Pull-backs of strata}
The pull-backs 
of  tautological classes to the boundary may now be
explicitly determined.  The basic calculation is 
the pull-back of the fundamental class of one boundary stratum
to another.  In terms of the diagram of Section \ref{specl}, we want
to compute $\xi_A^*(\xi_{B*}[\M_B])$.  As
we have identified $F_{A,B}$ explicitly as a smooth stack, the
pull-back will be straightforward to compute.  The 
intersection product is a sum of contributions of each component
of $F_{A,B}$, and
each  contribution is the Euler class of an excess bundle
on the component.  

The components of $F_{A,B}$ have been
identified in Proposition \ref{freddy}. Let
$C$ be a generic $(A,B)$-graph, and let 
$\M_C$ be the corresponding component of $F_{A,B}$. 
The excess bundle is easily identified on $\M_C$.
First, we observe that the normal bundle to $\xi_A$ naturally
splits as a copy of $e(A)$ line bundles. Let the edge $e$ be the join of
the distinct half-edges $h,h'$ incident to the vertices
$v,v'$ (which  may coincide).
The line bundle associated to $e$ is 
$$T_h \otimes T_{h'}$$ 
where $T_h$ and $T_{h'}$ are the tangent lines at
$h$ and $h'$ of the
factors 
$\M_{g(v),n(v)}$ and $\M_{g(v),n(v')}$ respectively.  
The normal bundle to $\M_C$ in $\M_A$ is a sum
of the analogous line bundles for those edges of $C$ which do not correspond
to edges of $A$.  
Precisely the same situation holds with respect to
$B$.  We can conclude that the excess
normal bundle of $\M_C$ thought of as a component of $F_{A,B}$ is
exactly the sum of the line bundles corresponding to those edges of
$C$ which correspond to edges of both $A$ and $B$.

We have deduced 
the following formula: 
\begin{equation}
\label{ffxx}
\xi_A^*(\xi_{B*}([\M_B])=\sum_C 
\xi_{C,A*}(\prod_{e=h+h'} -\pi_{v}^*(\psi_h) - \pi_{v'}^*(\psi_{h'})).
\end{equation}
The sum is over all generic $(A,B)$-graphs $C$.  The product is over
all edges $e$ of $C$ which come from both an edge of $A$ and an edge of $B$, and
$v,v'$ are the vertices joined by $e$. The morphism $\xi_{C,A}$ denotes
the natural map from $\M_C$ to $\M_A$.

Formula (\ref{ffxx}) yields an explicit 
tautological K\"unneth decomposition of the pull-back 
class since the morphism $\xi_{C,A}$ is simply the product of various
boundary strata maps over the factors of $\M_A$.  

We will compute a simple example to illustrate the formula.  Consider
the boundary divisor $\Delta_0$ in $\M_g$ corresponding to
the morphism $$i:\M_{g-1,2}\rarr \M_g.$$  The graph $A$ of
$\Delta_0$ has one vertex of genus 
$g-1$ and one self-edge.  
We will compute the
self intersection of the stratum: 
$\xi_A^*(\xi_{A*}[\M_A])$.  

We first write down all generic $(A,A)$-graphs.  
They are $A$ itself with the obvious $(A,A)$-structure, and then one
graph for each integer from $0$ to 
$\lfloor \frac{g-1}{2}\rfloor$.  $C_0$ is the
graph with one vertex of genus $g-2$ and two loops.  $C_0$ has two distinct
isomorphism classes of $(A,A)$-structures, but 
only one of them is generic: the $(A,A)$-structure
where the edge contracted for the first $A$-structure is different
from the edge contracted for the second $A$-structure.  Similarly, $C_i$ is
the graph with a vertex of  genus $i$ and another vertex of genus $g-i-1$ 
connected
to each other by two edges.  The unique generic $(A,A)$-structure is obtained
by contracting a different edge for the two $A$-structures.  Applying
formula (\ref{ffxx}), we find:
$$\xi_A^*(\xi_{A*}([\M_A]))= -\psi_1-\psi_2 + \xi_{0*}([\M_{g-2,4}]) + 
\sum_{i=1}^{\lfloor\frac{g-1}{2}\rfloor} \xi_{i*}([\M_{i,2} \times \M_{g-i-1,2}])
$$
with hopefully evident notations.  

Notice the boundary strata 
corresponding to $\M_{i,3}\times\M_{g-i-1,1}$ do not appear in the above formula
because the corresponding dual graphs do not admit a generic 
$(A,A)$-structure.  In more geometric terms, these strata do not contribute
an extra term because they have only one non-disconnecting node.

\subsection{Pull-backs of tautological classes}
We observe that our calculations easily generalize to computing
pull-backs of arbitrary tautological classes to boundary strata.

Define the tautological $\kappa$ classes by:
$$\pi_*(\psi_{n+1}^{l+1})=\kappa_l \in R^*(\overline{M}_{g,n}),$$
where $\pi$ is the map forgetting the last marking $n+1$. 
The first observation is the following result concerning the push-forwards
of the $\psi$ and $\kappa$ classes. 

\begin{pr}
\label{george}
Let $\pi: \overline{M}_{g,n+m} \rarr \overline{M}_{g,n}$ be the
map forgetting the last $m$ points. 
The $\pi$ push-forward of any element of the subring of
$R^*(\overline{M}_{g,n+m})$
generated by
$$\psi_1, \ldots, \psi_m, \{\kappa_i\}_{i\in {{\mathbb Z}_{\geq 0 }}}.$$
lies in the subring of $R^*(\overline{M}_{g,n})$ generated by
$$\psi_1, \ldots, \psi_n, \{\kappa_i\}_{i\in {{\mathbb Z}_{\geq 0 }}}.$$
\end{pr}
\noindent A proof  can be found in [AC].

We can now describe a set of additive generators for
$R^*(\overline{M}_{g,n})$. 
Let $B$ be a stable graph of genus $g$ with $n$ legs. 
For each vertex
$v$ of $B$, let
$$\theta_v\in R^*(\overline{M}_{g(v), n(v)})$$ be an
arbitrary monomial in the
cotangent line and $\kappa$ classes of the vertex moduli space.

\begin{pr}\label{maude}
 $R^*(\overline{M}_{g,n})$ is generated additively by classes of the form 
$$\xi_{B*}(\prod_{v\in V(B)} \theta_v).$$ 
\end{pr}
\bpf 
By the definition of $R^*(\overline{M}_{g,n})$, 
the claimed generators lie in the tautological ring.

We first show the span of the generators is closed under the intersection product.
The closure follows from:
\begin{enumerate}
\item[(i)] the pull-back formula (\ref{ffxx}) for strata classes,
\item[(ii)] the trivial pull-back formula for cotangent lines under
             boundary maps,
\item[(iii)] the pull-back formula for $\kappa$ classes under boundary maps [AC]
$$\xi_{B}^*(\kappa_i) = \sum_{v\in V(B)} \kappa_i.$$
\end{enumerate}

To prove the claimed generators span $R^*(\overline{M}_{g,n})$, we must
prove the system defined by the generators is closed under
push-forward by the forgetting maps and the gluing maps.

Closure under push-forward by the forgetting maps is a consequence of
Proposition \ref{george}. Closure under push-forward by the gluing maps
is a trivial condition.
\epf

\begin{cor}
$R^*(\overline{M}_{g,n})$ is a finite dimensional ${\mathbb Q}$-vector space.
\end{cor}
\bpf
The set of stable graphs $B$ for fixed $g$ and $n$ is finite, and
there are only finitely many non-vanishing monomials $\theta_v$ for
each vertex $v$.
\epf

\begin{pr}
Let $\gamma\in R^*(\overline{M}_{g,n})$. 
Let $A$ be a stable graph.
Let $$\xi_A: \overline{M}_A \rarr
\overline{M}_{g,n}.$$ Then, $\xi_A^*(\gamma)$ has a tautological K\"unneth
decomposition with respect to the product structure of $\overline{M}_A$.
\end{pr}
\bpf
The proposition follows from Proposition
\ref{maude} together with the pull-back formulas.
The pull-back formulas for the three types of classes all
yield tautological K\"unneth decompositions.
\epf

\vspace{+10 pt}
\noindent Department of Mathematics \\
\noindent Harvard University\\
\noindent Cambridge, MA 02138 \\
\noindent graber@@math.harvard.edu

\vspace{+10 pt}
\noindent
Department of Mathematics \\
\noindent California Institute of Technology \\
\noindent Pasadena, CA 91125 \\
\noindent rahulp@@cco.caltech.edu
\end{document}